\documentclass[12pt,letterpaper]{amsart}
\usepackage{palatino, euler, epic,eepic, amssymb, xypic, floatflt,
  microtype}
\usepackage{hyperref}
\input xy
\xyoption{all}

 \newlength{\baseunit}               
 \newcount{\numlines}                
\usepackage{pdfsync,url, mathrsfs}

\numberwithin{equation}{section}

\newcommand{\point}{\vspace{2mm}\par \noindent
  \refstepcounter{equation}{\theequation.} } 
\newcommand{\tpoint}[1]{\vspace{2mm}\par \noindent \refstepcounter{equation}{\theequation.}
  {\bf #1. ---} }
\newcommand{\epoint}[1]{\vspace{2mm}\par \noindent \refstepcounter{equation}{\theequation.}
  {\em #1.} }
\newcommand{\bpoint}[1]{\vspace{2mm}\par \noindent \refstepcounter{equation}{\theequation.}
  {\bf #1.} }
\newcommand{\ravispace}{ \vspace{2mm}}
\newcommand{\cut}[1]{}

\newcommand\hidden[1]{}

\newcommand{\bpf}{\noindent {\em Proof.  }}
\newcommand{\epf}{\qed \vspace{+10pt}}

\newcommand\al{\alpha}
\newcommand\be{\beta}

\newcommand\De{\Delta}
\newcommand\la{\lambda}
\newcommand\ka{\kappa}

\newcommand{\cK}{{\mathscr{K}}}
\renewcommand{\cL}{{\mathscr{L}}}
\newcommand{\proj}{\mathbb{P}}
\newcommand{\A}{\mathbb{A}}
\newcommand{\E}{\mathbb{E}}
\newcommand{\G}{\mathbb{G}}
\newcommand{\cB}{\mathscr{B}}
\newcommand{\cC}{\mathscr{C}}
\newcommand{\cE}{\mathscr{E}}
\newcommand{\cF}{\mathscr{F}}
\newcommand{\cG}{\mathscr{G}}
\renewcommand{\cH}{\mathscr{H}}
\newcommand{\cM}{\mathscr{M}}
\newcommand{\cm}{\mathscr{M}}
\newcommand{\cmbar}{\overline{\cm}}
\newcommand{\oh}{\mathscr{O}}

\newcommand{\cQ}{\mathscr{Q}}
\newcommand{\cT}{\mathscr{T}}
\newcommand{\cV}{\mathscr{V}}
\newcommand{\cW}{\mathscr{W}}
\newcommand{\Ext}{\operatorname{Ext}}
\newcommand{\Pic}{\operatorname{Pic}}
\newcommand{\rank}{\operatorname{rank}}
\newcommand{\Aut}{\operatorname{Aut}}
\newcommand{\Hom}{\operatorname{Hom}}

\newcommand{\Sym}{\operatorname{Sym}}

\newcommand{\Q}{\mathbb{Q}}

\newcommand{\Z}{\mathbb{Z}}
\newcommand{\F}{\mathbb{F}}
\newcommand{\C}{\mathbb{C}}

\renewcommand{\Q}{\mathbb{Q}}

\title{The Chow ring of the moduli space of curves of genus $6$}

\date{July 7, 2013.}
\author{Nikola Penev and Ravi Vakil}
\address{
N. Penev,  Dept.\ of Mathematics, Stanford University,
  Stanford, CA 94305 USA 
\newline \indent
R. Vakil,  Dept.\ of Mathematics, Stanford University,
  Stanford, CA 94305 USA  } 
\email{vakil@math.stanford.edu, npenev@gmail.com}

\begin{document}
\begin{abstract}  
  We determine the Chow ring (with $\Q$-coefficients) of $\cm_6$  by showing that all Chow classes are
  tautological.  In particular, all algebraic cohomology is
  tautological, and the natural map from Chow to cohomology is
  injective.  To demonstrate the utility of these methods, we also give quick
  derivations of the Chow groups of moduli spaces of curves of lower
  genus.  The genus $6$ case relies on the
  particularly beautiful Brill-Noether theory in this case, and in
  particular on a rank $5$ vector bundle ``relativizing'' a baby case
  of a celebrated construction of Mukai, which we interpret as a
  subbundle of the rank $6$ vector bundle of quadrics cutting out the
  canonical curve.
\end{abstract}
\maketitle
\setcounter{tocdepth}{1} 
\tableofcontents

\section{Introduction}

Modern progress in our algebro-geometric understanding  of  the moduli
space of smooth curves
was launched by  Faber's papers \cite{Fg3, Fg4} on the Chow rings
of $\cm_3$ and $\cm_4$, which described tools for 
understanding about Chow rings of $\cm_g$, and Faber's conjecture
\cite[Conj.~1]{Fconj}, which suggested the existence of a rich structure in the
``tautological'' part of the Chow ring $R^*(\cm_g) \subset
A^*(\cm_g)$.  
(Important note:  throughout this paper, Chow groups are taken with $\Q$-coefficients.)
  Earlier work of Mumford \cite{mumford},
and Witten's conjecture \cite{witten},
provided the foundations for these papers.  

Mumford
 \cite[p.~318]{mumford}
 earlier described the Chow ring of $\cm_2$,
and Izadi 
 \cite{izadi}
later determined the Chow ring of $\cm_5$.  In
genus up to $5$ the Chow ring is all tautological, and of the
form
\begin{equation}
A^*(\cm_g) = \Q[\ka_1]/ (\ka_1^{g-1}). \label{eq:g5}
\end{equation}
In this paper, we describe the Chow ring of $\cm_6$:

\tpoint{Main Theorem}{\em   
\begin{equation}
\label{eq:RM6}
A^*(\cm_6) = \Q[ \ka_1, \ka_2] / ( 127 \ka_1^3 - 2304 \ka_1 \ka_2, 113
\ka_1^4 - 36864 \ka_2^2).
\end{equation}
In particular, in $H^*(\cm_6, \Q)$, all ``algebraic cohomology'' is
tautological, and the natural map $A^*(\cm_6) \rightarrow
H^{2*}(\cm_6, \Q)$
is an injection.\label{t:main}}\ravispace

The explicit description \eqref{eq:RM6}  of the ring hides its elegance.  It should be
understood  as follows.  Faber's conjecture \cite{Fconj} describes one ring for each
genus $g$, determined by the fact that it is a Poincar\'e duality ring
of dimension $g-2$, with known generators, and relations determined by
the top intersections of these generators, which in turn have a
beautiful description in terms of the symmetric group.  This ring
agrees with the tautological ring in genus up to $g=23$.  The case
under consideration here is the first genus in which the ring is not
of the form \eqref{eq:g5} (see \cite[p.~123]{Fconj}).

Theorem~\ref{t:main} is proved in \S \ref{s:background6} and \S \ref{s:mukai}; we make
the argument explicit in \S \ref{pf:main}.  

\bpoint{New ingredients}
We briefly describe the new points of view which make
Theorem~\ref{t:main} possible.

{\em (i)} We  take advantage of the fact that Faber has already determined
the tautological ring \cite[p.~123]{Fconj}; we will show that all Chow
classes are
tautological, and need not worry about determining relations.

{\em (ii)} We cut up $\cm_6$ into locally closed strata as is
traditional, but we do not bother computing the Chow groups of the
strata.  Instead, we choose strata which are group quotients, and use
a versatile theorem of Vistoli (Theorem~\ref{t:vistoli}) to show that the Chow
rings are generated by Chern classes of some natural vector bundle(s).
We then show that the fundamental class of the stratum is
tautological, and also relate the Chern classes of the vector bundle(s)
to the Hodge bundle to show that they too are restrictions of
tautological classes on $\cm_6$.

{\em (iii)}  For a large open subset $\cm_6^M$ of $\cm_6$ (the ``Mukai
locus'', the complement of the trigonal, hyperelliptic, and plane
quintic loci), we use Mukai's work
describing the corresponding canonical curves as complete intersections in
$G(2,5)$, and in particular we construct a corresponding rank $5$ vector bundle $\cV$
on this open subset.  Several aspects of Mukai's construction come into play.

The fact that it is tractable to determine $A^*(\cm_6)$ requires a 
number of fortunate coincidences.  But we hope that 
our methods will be useful in other circumstances.  As a first
example, it seems plausible that such methods can show that
$A^*(\cm_g)$ is finitely generated for $7 \leq g \leq 9$, using Mukai's
description of large open subsets of $\cm_g$ in this genus range.
(There seems no reason to believe that $A^*(\cm_g)$ is finitely
generated in general.)

\epoint{Example} For those less familiar with the tautological ring,
we point out as an example the unusual way in which classes are shown to be zero.  We explain how we know that the fundamental class
of the locus $\cB_6$ of bi-elliptic curves (those curves admitting a
double cover of an elliptic curve) is $0$.  
(The difficulty in showing this fact obstructed an earlier attempt by
the authors to prove the Main Theorem~\ref{t:main}.)
Even the fact that $[\cB_6]$ is
$0$ in cohomology is not clear; this is the first in a series of
interesting questions on the cohomology of $\mathscr{A}_g$
(S. Grushevsky, private communication).  

The
vanishing of $[\cB_6]$  is for the following  unusual reasons (which are  most unlikely to
generalize).

The locus $\cB_6$ is irreducible of codimension $5$.  By Looijenga's
theorem (the main result of \cite{Looijenga-Inventiones}, see  \S
\ref{s:taurin}(i)), $R^5(\cm_6)=0$, so it suffices to show that
$[\cB_6]$ is tautological.  This class lies in the open subset
$\cm_6^M$ of Mukai-general curves (see
Theorem~\ref{t:frommukai}).  We describe $\cm_6^M$ as a quotient by
$GL(5)$ (Theorem~\ref{t:mukai}).  The dividend is a hypersurface
complement in an affine bundle over a Grassmannian.  The Chow groups
of the quotient are generated, in a suitable sense, by (i) the Chow
groups of the Grassmannian, and (ii) the Chern classes of the vector
bundle corresponding to the $GL(5)$-quotient (Prop.~\ref{p:fv}).  These can be related to
each other, and the Hodge bundle, to show they are all tautological.
Thus $[\cB_6]$ is forced to be tautological, and hence $0$.

\bpoint{Structure of paper} In \S \ref{s:background}, we collect the
tools we will use.  In \S \ref{s:g2345}, we demonstrate the approach,
by quickly recovering the Chow rings of $\cm_2$ through $\cm_5$.  In
\S \ref{s:background6}, we give background on $\cm_6$, mentioning in
particular the facts we need.  In \S \ref{s:trig}, we show that all
classes supported on the trigonal, hyperelliptic, or plane quintic
locus of $\cm_6$ are tautological.  In \S \ref{s:mukaigeneral}, we describe the rank
$5$ ``Mukai'' bundle on the Mukai locus $\cm_6^M$.  Finally, in \S \ref{s:mukai}, we show
that all Chow classes supported on the Mukai locus are tautological.
This concludes the proof of Theorem~\ref{t:main}

\bpoint{Acknowledgments} This paper owes an obvious debt to C. Faber,
and we are grateful to him.  We also thank G. Farkas and S. Grushevsky
for helpful conversations.  The second author was supported by NSF
grants DMS-1100771 and DMS-1159156, and the Simons Foundation.  He
thanks the American Institute of Mathematics for hospitality during
the writing of this article.

\section{Background}
\label{s:background}

We work over $\C$ for convenience.  All quotients should be
interpreted as stacks,
not coarse moduli spaces.
If $V$ is a vector space, then $\proj V$ is the
projectivization in the sense of Grothendieck, i.e., 
the  space of
one-dimensional {\em quotients} of $V$.  Thus $V$ is the space of
``linear functionals'' on $\proj V$.
 We work throughout with finite type
schemes or Deligne-Mumford stacks.
The tautological ring $R^*(\cm_g)$ is the subring of $A^*(\cm_g)$ (not,
a priori, $H^*(\cm_g, \Q)$).  If $U$ is an open substack of $\cm_g$,
then $R^*(U)$ is defined as $R^*(\cm_g)|_U$. 

\bpoint{Background on Chow groups}

\label{s:chow} We say that an irreducible variety $X$ is {\em
  Chow-free} if $A_* X$ is generated (as a $\Q$-vector space) by the
fundamental class $[X]$.  If $Z \hookrightarrow X$ is a closed
embedding with complement $U$, then we have an exact sequence of Chow
groups 
\begin{equation}\label{eq:excision}A_* Z \rightarrow A_* X \rightarrow
  A_* U \rightarrow 0\end{equation}
(see \cite[Prop.~1.8]{F} for schemes, and \cite[Prop.~2.3.6]{kresch}
for the extension to stacks).
Thus if $D \subset \proj^n$ is a closed
subset containing a (nonempty) hypersurface, then $\proj^n \setminus
D$ is Chow-free.  Also, to show that Chow groups of $\cm_g$ are
tautological, it suffices to stratify $\cm_g$, and to inductively show that the
pushforwards of Chow groups of the closures of these strata are all tautological.

For the remainder of \S \ref{s:chow}, suppose $\pi: X \rightarrow Y$ is
a morphism.

\point \label{s:chow.I} (i) If $\pi$ is finite and flat of degree $d$,
then $\pi^*: A_* Y \hookrightarrow A_* X$ is an injection, as $\frac 1
d \pi_*$ is a one-sided inverse.

(ii)  It is straightforward to show that if $\pi$ is a gerbe banded
by a finite group, then $\pi^*$ is an isomorphism.  (The cycles on $X$
and on $Y$ are obviously identified.)

\tpoint{Theorem (Vistoli)}
\label{t:vistoli}
{\em Suppose $G=GL$
or $SL$, and $\pi$ is a quotient by $G$ (i.e., $Y = X/G$
categorically).  Let $V$ be the corresponding vector bundle $V$ on
$Y$; if $G= SL$, then $c_1(V)=0$.  Then: 
\begin{enumerate}
\item[(i)] The pullback $\pi^*: A_* Y
\twoheadrightarrow A_* X$ is surjective.  
\item[(ii)] The kernel of $\pi^*$ is
generated as a group by classes of the form $c_i(V) \cap Z$,  where
$i>0$, and $Z \in  A_*(Y)$.
\item[(iii)] 
If $X$ (and hence $Y$) is smooth, then $A^*Y$ is generated as a ring
by $c_1(V)$, \dots, $c_n(V)$ and (any chosen) lifts of
generators of $A^* X$.  In particular, 
if $X$ is Chow-free,
then $A_* Y$ is generated by $\{ c_i(V) \}$.
\end{enumerate}}

Parts (i) and (ii) are 
\cite[Thm.~2]{V}, and (iii) follows by induction on codimension.

\bpoint{Background on the tautological ring}

\label{s:taurin}(i) Let $\cH_g \subset \cm_g$ be the locus of hyperelliptic curves.
We will use the following, from the main result of
\cite{Looijenga-Inventiones}: $R^i(\cm_g) = 0$ for $i>g-2$, and
$R^{g-2}(\cm_g)$ is generated by $[ \cH_g ]$.   (In fact, $[\cH_g]
\neq 0$, \cite[l.~3]{Fnonvan},  but we will not need
this.)

(ii)  If a Brill-Noether locus is of the expected dimension, then its
fundamental class is tautological.
 (The cases of interest to us here are the loci of
hyperelliptic curves, trigonal curves, and plane quintics.)
This was shown in \cite[p.~116-118]{Fconj}, by suitably describing each
Brill-Noether locus as a
degeneracy locus. 

(iii) 
{\em $\Pic \cm_g \otimes_{\Z} \Q =  \Q[\la]$.} \label{picmg}
This follows from the section determining $\Pic \cmbar_g$ in
\cite{M2}. 
(In the genera of interest to us, $g \leq 6$, purely algebraic
arguments are possible.  For example, in the case $g=6$, an
argument is given in the proof of Proposition~\ref{p:picm6}.)

\section{Moduli of curves of small Clifford index}

\label{s:g2345}In this section, we show that Chow groups of various classes pushed
forward from moduli spaces of curves of low Clifford index are tautological.  We first
discuss hyperelliptic and trigonal curves.  We then deal with low
genus curves, partially to show the utility of the methods.  The
trigonal case is the main substantive part of this section.

\bpoint{Hyperelliptic curves}

\label{s:hyp}(This discussion can be done simply in other ways, but we
deliberately use the same methods as we use later.)  Fix a genus $g \geq 2$.  
Recall that $[\cH_g]$ is tautological, by  \S
\ref{s:taurin}(i) (Looijenga's Theorem) or \S \ref{s:taurin}(ii).
Thus to show that the pushforward of the Chow groups of $\cH_g$ to
$\cm_g$ are tautological, it suffices to show that the Chow ring of
$\cH_g$ is generated by elements that are pulled back from
elements of $R^*(\cm_g)$. 

Let $\cH'_g$ be the moduli space of $2g+2$ unordered distinct points
on a genus $0$ smooth curve.  Then there is a natural map $\rho: \cH_g
\rightarrow \cH'_g$, which is a $\mu_2$-gerbe.  The
pullback map $\rho^*:  A^*(\cH'_g) \rightarrow A^*(\cH_g)$ is an isomorphism, by \S \ref{s:chow.I}(ii).  Thus it
suffices to show that the Chow ring of $A^*(\cH'_g)$ is tautological,
i.e., its elements pull back to tautological elements of $\cH_g$.
Now $\cH'_g = (\C^{2g+3}
\setminus \Delta ) / GL(2)$, where $\C^{2g+3}$ parametrizes degree
$2g+3$ homogeneous polynomials in $\C[x,y]$, and $\De$ is the discriminant locus.    The
$GL(2)$-quotient induces a vector bundle $\cW$ on $\cH'_g$.  By
Vistoli's Theorem~\ref{t:vistoli}(iii), as $\C^{2g+3} \setminus \De$ is
Chow-free, $A^*(\cH'_g)$ is generated by $c_1(\cW)$ and
$c_2(\cW)$.  Now $\rho^* (\proj \Sym^{g-1} \cW)$ is
identified with the target of the canonical embedding $\proj \E$
(where $\E$ is the Hodge bundle), so $\Sym^{g-1} \cW$ is
isomorphic to $\E \otimes \cL$ for some invertible sheaf $\cL$ on
$\cH_g$.  But $\Pic_{\Q} \cH_g = 0$.   (This is well known, and can
even be shown by the methods described here; see
\cite[Thm.~5.1]{av} for a much more general statement.) 
 Thus $\Sym^{g-1} \cW \cong \E$, so 
$\Sym^{g-1} \cW$ has
tautological Chern classes in $A^*(\cH_g)$, from which (by considering
Chern roots) $\cW$ does as well.    

Thus the pushforward of the Chow groups of $\cH_g$ are all
tautological.
(In particular, by \S \ref{s:taurin}(i), the pushforward yields just
multiples of the fundamental class $[\cH_g]$ in $A^*(\cm_g)$.)

\bpoint{Trigonal curves}

\label{s:trigonal}It would be good to know that the Chow groups of 
the locus of trigonal curves are tautological in $\cm_g$.
We will establish something weaker.  See \S \ref{s:tritau} for what is
missing.

Fix a genus $g \geq 4$.  Let $\cT_g$ be the stack of smooth curves,
along with the data of a triple covers of a genus $0$ curve.  If $g>4$, the natural moduli morphism
$\cT_g \rightarrow \cm_g$ is a closed embedding.  (This is
well known and straightforward to show;
simply
describe an equivalence of categories between $\cT_g$ and a closed
substack of $\cm_g$.  See \cite[Prop.~2.2]{bv} for details, and much
more.)

Recall that $\cT_g$ is
stratified by {\em Maroni invariant}.  The Maroni invariant of a
trigonal cover $\pi: C \rightarrow \proj^1$ is the unique nonnegative
integer $n$ such that  $\pi$ factors through a closed embedding $\pi:
C \hookrightarrow \F_n$ where $\F_n$ is the Hirzebruch surface
$\proj_{\proj^1} ( \oh \oplus \oh(n) )$ (with the implicit projection
to $\proj^1$ as the trigonal structure morphism).  We denote the ``standard'' elements of $A^1(\F_n)$ by
$E$ (the class of a section of self-intersection $-n$, unique if
$n>0$), $S$ (the class of a
section not meeting $E$, of self-intersection $n$), and $F$ (the fiber
class). 
Important note:  we will use $S$ to mean an actual choice of section, not
just a Chow class.

 Then $C \hookrightarrow \F_n$ is in class $3S + kF$, where $k
= (g-3n+2)/2$.  In particular, $n$ must have the same parity as $g$,
and $n \leq (g+2)/3$.

\tpoint{Theorem} {\em    If $\cT_{g,n}$ is the space  of trigonal  genus $g$ smooth
curves of Maroni invariant $n$, then any Chow class on $\cT_{g,n}$ is
the restriction of a tautological class on $\cm_g$.}\label{t:trigonal}

\point \label{s:tritau}If we knew that the fundamental classes
$[\cT_{g,n}]$  were {\em also} tautological
on $\cm_g$ for all $n$, we would then know that the (pushforward of the) Chow
groups of the locus of {\em all} trigonal curves are tautological in
$\cm_g$.

\epoint{Proof}   Consider
$\proj^N := \proj H^0(\F_n, \oh(3S+kF))$, where $k= (g-3n+2)/2$.  
Let $\De$ be the divisor on $\proj^N$ corresponding to  singular
curves.  (Note
that $\proj^N \setminus  \De$ is
Chow-free.)  Then $\cT_{g,n} \cong (\proj^N \setminus \De) / \Aut
\F_n$. 
Our goal is to show that the Chow groups of this
quotient are all restrictions of tautological classes.

Let $G \subset \Aut \F_n$ be the subgroup 
fixing the curve $S$
as a set (but not necessarily pointwise).  Then 
$(\proj^N \setminus \De) / G \rightarrow \cT_{g,n}$ is an
$\A^{n+2}$-bundle (as $h^0(\F_n, \oh(S)) = n+2$), so pullback by this morphism induces an isomorphism
of Chow rings
by \cite[Cor.~2.5.7]{kresch}
 Hence it suffices to show that the Chow
groups of 
$(\proj^N \setminus \De) / G$
 are restrictions of tautological classes.

We have a
filtration
$$
1 \rightarrow \G_m \rightarrow G \rightarrow \Aut (\proj^1) \rightarrow 1
$$
where the $\G_m$ acts on $\F_n$ by its action on the fibers (fixing
both $E$ and $S$), and $\Aut(\proj^1)$ is the automorphism group  of the target
$\proj^1$.    
We extend $G$ to a group $G'$ along the surjection $SL(2) \rightarrow
\Aut(\proj^1) = PGL(2)$:
$$
1 \rightarrow \G_m \rightarrow G' \rightarrow SL(2) \rightarrow 1.
$$
Then 
$(\proj^N \setminus \De) / G' \rightarrow  (\proj^N \setminus \De) /
G$
is a $\mu_2$-gerbe, so pullback by this morphism induces
an isomorphism of Chow groups by  \S \ref{s:chow.I}(ii).
Hence it suffices to show that the Chow
groups of 
$(\proj^N \setminus \De) / G'$
 are restrictions of tautological classes.

By Vistoli's Theorem~\ref{t:vistoli}(iii), the Chow groups of this quotient $(\proj^N
\setminus \De) / G'$ will be generated by $\al_1$, which is  defined
to be the first
Chern class of the line bundle $\cL$ corresponding to $\G_m$, and
$\be_2$, which is defined to be the second Chern class of the vector bundle $V$
corresponding to the $SL(2)$.  (Here we use $c_1(V)=0$.)
In particular,  the only line bundles on the quotient $(\proj^N \setminus
\De) / G'$ are $\Q$-multiples of $\al_1$.

Geometrically, $V$ should be
interpreted as a vector space such that $\proj V$ is the target
$\proj^1$; and $\cL$ can be informally interpreted as a fiber of $\F_n
\setminus E$ over any point of the target $\proj^1$ (with the
intersection with $S$ as $0$).

We now relate these two bundles  to the Hodge bundle $\E$.  We motivate this with
some geometry.  If we canonically embed the trigonal curve in
$\proj^{g-1}$, the ambient Hirzebruch surface ``comes along for the
ride'' as a scroll.  Geometrically:  the fibers of the $\proj^1$ are recovered
from the fact that the ``trigonally conjugate points'' are collinear
in the canonical embedding.
Sheaf-theoretically, this is the fact that if $D$ is such a trigonal
curve on $\F_n$, then 
$$H^0(\F_n, \oh(S + (n+k-2) F)) \rightarrow
H^0(D, \oh_D(S+(n+k-2) F))$$ 
is an isomorphism, and the right side are canonically
the differentials on $D$.

Better yet, since we have a privileged curve $S$ on $\F_n$ (not just a
curve class), we have a canonical direct sum decomposition
$$H^0(\F_n, \oh(S + (n+k-2) F)) 
 = H^0(S, \oh(2n+k-2) ) \oplus
H^0(E, \oh(n+k-2) )$$ which we interpret geometrically as
follows.  We have $\F_n$ embedded in $\proj^{g-1}$, and $E$ is a
degree $n+k-2$ rational normal curve, and $S$ is a degree $2n+k-2$
rational normal curve, and the projective spaces they span have empty
intersection, and together span the full $\proj^{g-1}$.

Let $\proj^1 = \proj V$, where $V$ is a two-dimensional vector space,
and let $\cL$ be the one-dimensional tautological representation of $\G_m$.
Then the canonical embedding is given by 
\begin{equation}
\label{eq:this}\left( \Sym^{2n+k-2} V \right) \oplus \left( \cL \otimes \Sym^{n+k-2} V \right).\end{equation}
Here $\be_2 = c_2(V)$, and $\al_1 = c_1(\cL)$.
Since we have identified  projectivization of this bundle
\eqref{eq:this}  with $\proj \E$, 
the Hodge bundle $\E$ must differ from this bundle \eqref{eq:this} by some line bundle $\cM$:
$$\E \otimes \cM \cong 
\left( \Sym^{2n+k-2} V \right) \oplus \left( \cL \otimes \Sym^{n+k-2} V \right).$$
As mentioned above, $\cM$ must be some $\Q$-multiple of $\cL$ (i.e., 
$c_1(\cM)$ is some multiple of $\al_1$).
Comparing the first Chern class of both sides, we see that $\cL$  is
some multiple of $\la_1$. Comparing the second Chern class of both
sides, we see that $\be_2$ is a linear combination of $\la_1^2$ and $\la_2$. \epf

\bpoint{Curves of low genus}

We now apply these methods to show  that $A^*(\cm_g)$ is tautological
for $g<6$.
The case $g=2$ is a special case of our hyperelliptic discussion, \S
\ref{s:hyp}.

\epoint{Genus $3$}

We stratify $\cm_3$ into the nonhyperelliptic locus $\cm_3 \setminus \cH_3$, and the
hyperelliptic locus $\cH_3$.  We established in \S \ref{s:hyp}
that all Chow groups of $\cH_3$ are tautological (and indeed that
$\cH_3$ is Chow-free), so we concentrate on the open stratum
$\cm_3 \setminus \cH_3$, which parametrizes smooth plane quartics up to
automorphisms of the plane.  Hence $\cm_3 \setminus \cH_3 \cong (H^0(\proj^{2},
\oh(4) )\setminus \Delta) / GL(3)$.  The vector bundle corresponding
to the $GL(3)$-quotient is
the Hodge bundle.  Thus the Chow ring of $\cm_3 \setminus \cH_3$ is generated by the
Chern classes of the Hodge bundle, so $A^*(\cm_3 \setminus \cH_3)$ is
all tautological.

\epoint{Genus $4$}

\label{s:g4}We cut  $\cm_4$ into three strata.  

\noindent {\em $\cm_4^0$ (dense, dimension $9$):  those curves whose canonical models are the
intersections of smooth quadrics and cubics.}  
This has a double   cover 
(i.e., degree $2$ finite flat cover)
$\cm \rightarrow \cm_4^0$ corresponding to
such curves with a choice of one of the two $g^1_3$'s.    Then by our
trigonal discussion \S \ref{s:trigonal}, $A^*(\cm)$ is generated by tautological classes.
Hence by \S \ref{s:chow.I}(i), $A^* (\cm_4^0)$ is all tautological.

\noindent {\em$\cT_{g,2}$ (dimension $8$): those curves whose canonical models are
  intersections of cones and cubics (or:  the divisor of curves with a
vanishing theta-null).}
We first note that $[\cT_{g,2}] \in A^1(\cm_4)$ is tautological, by \S
\ref{s:taurin}(iii).
(As noted there, this can be shown in many purely algebraic ways.)
   We are then done by Theorem~\ref{t:trigonal}.

\noindent {\em $\cH_4$:  The hyperelliptic locus (dimension $7$).}  This was dealt with in
\S \ref{s:hyp}. 

\epoint{Genus $5$}

We cut $\cm_5$ into  three strata.\label{s:g5} 

\noindent {\em  $\cm_5^0$ (dense, dimension $12$): non-trigonal
  non-hyperelliptic curves.}
These are the curves whose canonical models are complete intersections
of three quadrics in $\proj^4$.  We consider $G(3,15)$ parametrizing
nets (linearly embedded $\proj^2$'s) of quadrics in $\proj^4$.  This
comes with a rank $3$ bundle, the tautological bundle $\cQ$ of the
Grassmannian.  
Let $\De$  be the divisor on $G(3, 15)$ corresponding to nets of quadrics
whose intersection is not a smooth curve.  
Then $(G(3,15) \setminus \De) / PGL(5) \cong \cm^0_5$ (we have a
equivalence of moduli problems).
We then consider $(G(3,15) \setminus \De)/SL(5)$.   This is a
$\mu_5$-gerbe over $(G(3,15) \setminus \De) / PGL(5) \cong \cm^0_5$,
and thus has the same rational Chow groups by \S \ref{s:chow.I}(ii).

Let $\cV$ be the tautological rank $5$ vector bundle on the quotient
corresponding to the $SL(5)$-bundle.  By Vistoli's Theorem~\ref{t:vistoli}, the Chow classes of the quotient are generated
by (lifts of) the cohomology classes on $G(3,15)$, and the Chern
classes of $\cV$.  Recall that the Chow ring of $G(k,n)$ is generated
by the Chern classes of the tautological bundle (see for example
\cite[\S 14.7]{F}).

Now $\proj \E \cong \proj \cV$ as projective bundles over $\cm_5^0$,
as the universal curve is embedded canonically in both projective
bundles.  Thus $\cV\otimes \cL \cong \E$, for some line bundle $\cL$
on $\cm_5^0$.  As $\det \cV = \oh$, $\cL$ is a (canonical) fifth
root of $\det \E$.  (Hence $c_1(\cL) = \la_1/5$ in $A^1 ((G(3,15)
\setminus \De)/SL(5))$.)  Now $\cQ$ is a rank $3$ subbundle of $\Sym^2
\cV$, so $\cQ \otimes \cL^{\otimes 2}$ is a subbundle of $\Sym^2 \E$.
This bundle $\cQ \otimes \cL^{\otimes 2}$ is the kernel of the
surjective map of vector bundles $\Sym^2 \pi_* \omega \rightarrow
\pi_* \omega^{\otimes 2}$, where $\pi: \cC_5^0 \rightarrow \cm_5^0$ is
the universal curve, and $\omega$ is the relative dualizing sheaf.
The Chern classes of $\Sym^2 \pi_* \omega = \Sym^2 \E$ are 
tautological, and the Chern classes of $\pi_* \omega^{\otimes
  2}$ 
(by Grothendieck-Riemann-Roch, see for example
\cite[p.~111]{Fconj}).
Thus the Chern classes of $\cQ$ are also tautological.

In conclusion,  the Chow
classes of the quotient are thus tautological, because they arise from
(lifts of) the cohomology classes on $G(3,15)$ (which are all
tautological), and the Chern classes of the $\cV \cong \E \otimes
\cL^\vee$.

\noindent {\em $\cT_{3,1}$ (dimension $11$): trigonal curves.}  All
genus $5$ trigonal curves have Maroni invariant $1$.  
As $[\cT_{3,1}] \in A^*(\cm_5)$ is tautological by \S
\ref{s:taurin}(ii), we are done by \S \ref{s:tritau}.

\noindent {\em $\cH_5$ (dimension $9$)  hyperelliptic curves.}
  This was dealt with in
\S \ref{s:hyp}.

\section{Stratification of $\cm_6$, and structure of the argument}

\label{s:background6}The following stratification follows immediately
from \cite[Thm.~5.2]{acgh} (Mumford's refinement of Marten's Theorem),
and straightforward dimension counts.

\tpoint{Theorem}
{\em \label{t:6stra}The space $\cm_6$ may be partitioned into the following locally closed
subsets.
\begin{enumerate}
\item [(i)]  {\em (dimension 15)} The locus of ``Brill-Noether-general'' curves, which we denote $\cm_6^{BN}$,
consisting of curves with finitely many $g^1_4$'s.
\item[(ii)]  {\em (dimension 13)} The  trigonal locus $\cT_6$.
\item[(iii)] {\em (dimension 12)}  The locus $\cQ_6$ of 
  plane quintics.
\item[(iv)] {\em (dimension 11)}  The hyperelliptic locus $\cH_6$.
\item[(v)] {\em (dimension 10)} The bi-elliptic locus $\cB_6$.
\end{enumerate}}

\tpoint{Proposition}
{\em  The union $\cm_6^M := \cm_6^{BN} \cup \cB_6$ is an open subset of $\cm_6$.
Equivalently, $\cB_6$ does not meet the closure (in $\cm_6$) of $\cT_6$, $\cQ_6$,
or $\cH_6$.}\ravispace

We call $\cm_6^M$ the Mukai locus (for reasons that will be clear by the
start of \S \ref{s:mukaigeneral}).

\bpf From the main theorem of \cite{schreyer} (see \cite[Table
1]{schreyer}), there are three possible Betti-tables for the canonical
rings smooth genus $6$ curves, depending on whether the curve is in
$\cm^M_6 = \cm^{BN}_6 \cup \cB_6$, $\cT_6 \cup \cQ_6$, or $\cH_6$.  By
upper-semicontinuity, from this table the first locus is open.

(Alternatively, one could consider $W^1_4(C)$ as $[C]$ varies over
$\cm_6$.  Or one could use the Enriques-Petri Theorem, see \cite[Main
Thm.]{saint-donat}.)  \epf

The following fact is a special case of \S \ref{s:taurin}(iii).

\tpoint{Proposition}
{\em $\Pic \cm_6 \otimes_{\Z} \Q =  \Q[\la]$.} \label{p:picm6}\ravispace

As promised there, we give a purely algebraic argument.

\bpf Let $\cm_6^0$ be the open subset of $\cm_6^{BN}$ corresponding to
those curves with precisely $5$ $g^1_4$'s.  It has a degree $120$
cover corresponding to such curves with $5$ {\em labeled} $g^1_4$'s.
These in turn are parametrized by a hypersurface complement in
$\proj^{15}$, parametrizing plane sextic curves, smooth except for
nodes at $[1,0,0]$, $[0,1,0]$, $[0,0,1]$, and $[1,1,1]$.  This has
trivial codimension $1$ rational Chow group, and thus $\cm_6^0$ does
too, by \ref{s:chow.I}(i).  The complement is an irreducible
divisor, generically corresponding to a genus $6$ curve with precisely
four $g^1_4$'s (one with multiplicity $2$); see
\cite[Exercises~V.A]{acgh}.  Finally, $\la$ is nonzero, for example by
using ampleness, or the fact that the tautological ring is known
\cite[p.~123]{Fconj}.  \epf

\bpoint{How we will prove the Main Theorem~\ref{t:main}}
\label{pf:main}As with the examples of $\cm_2$ through $\cm_5$ in \S
\ref{s:g2345}, we will prove the Main Theorem~\ref{t:main} by proving
that all Chow groups are tautological, and in turn do this by
showing that the pushforward of all Chow classes from all strata
described in Theorem~\ref{t:6stra} are tautological.  We deal with the
complement of the Mukai locus in the remainder of \S \ref{s:background6}, and we
deal with the Mukai locus in the remaining two sections of the paper.
By  the ``excision exact sequence for Chow groups'' \eqref{eq:excision}, this
will complete the proof of Theorem~\ref{t:main}.

\bpoint{Classes pushed forward from the complement of the Mukai locus
  are tautological} 
\label{s:trig}

We established in \S \ref{s:hyp} that the (pushforward of) Chow
classes on the hyperelliptic locus $\cH_6$ are tautological.

\epoint{The plane quintics $\cQ_6$}
by \S
\ref{s:taurin}(ii), $[\cQ_6]$ is tautological, 
 As $h^0(\proj^2, \oh(5)) = 21$, we may present
$\cQ_6$ as $(\A^{21} - \Delta) / GL(3)$, where $\De$ is the
``discriminant'' divisor
corresponding
to singular quintics.
Thus the Chow ring of $\cQ_6$ is generated by the Chern classes of the corresponding rank $3$
vector bundle $V$.  But $\Sym^2 V$ is identified with the Hodge bundle $\E|_{\cQ_6}$
(the canonical bundle of a plane quintic $C$ is canonically identified with
the restriction of $\oh_{\proj^2}(2)|_C$), so we can write
the Chern classes of $V$  in terms of the Chern classes $\la_i$ of the
Hodge bundle.

\epoint{Trigonal curves $\cT_6$}
Similar to the genus $4$ case (\S \ref{s:g4}), we further stratify $\cT_6$ by  Maroni-invariant:
it is the union of the Maroni-general locus $\cT_{6,0}$ and the
``Maroni divisor'' $\cT_{6,2}$.  
The first has tautological fundamental class by \S \ref{s:taurin}(i),
and the second does by \cite[equ.~(1.4)]{stankova}.
Then we are done by Theorem~\ref{t:trigonal}.

\section{The rank $5$ Mukai bundle on $\cm_6^M$}

\label{s:mukaigeneral}Mukai gives a remarkable characterization of ``most'' canonical
curves of genus $6$ through $9$ in his series of papers starting with
\cite{mukai}.   We will rely on his genus $6$
results, given in \cite[\S 5]{mukai}.   

\tpoint{Theorem (Mukai, \cite[\S 5]{mukai})} {\em If $C$ is a
  Mukai-general curve of genus $6$, then there is a unique stable rank
  $2$ vector bundle $E$ on $C$ with $\det E \cong \cK_C$ and having a
  $5$-dimensional vector space $V$  of sections.  The bundle $E$ is
  globally generated, and thus yields a morphism $C \rightarrow
  G(2,V)$, which is a closed embedding.  This closed embedding is a
  complete intersection of four linear forms and one quadratic form,
  under the Pl\"ucker embedding $G(2,V) \hookrightarrow \proj
  (\wedge^2 V)$.}  \label{t:frommukai}\ravispace

(Caution: when reading \cite[\S 5]{mukai}, note that ``del Pezzo
surfaces'' there are not assumed to be smooth.  More precisely, if $C$
is a Mukai-general curve, then it may be described as a quadratic
section of a dimension $2$ linear complete intersection of $G(2,5)$,
and this surface is a del Pezzo surface, possibly singular.  For a
complete description of the possible surfaces, see \cite[\S 5]{ah}.
For example, in the case of bi-elliptic $C$,  it is  a cone over an elliptic quintic, see
\cite[p.~169]{ah}.)

\point 
 \label{s:frommukai}As 
$\cK_{G(2,V)} = \oh_{\proj \wedge^2 V}(-5)|_{G(2,V)}$, by the adjunction formula, $\cK_C = \oh_{\proj \wedge^2
  V}(1)|_C$.    Let $F$ be the  four-dimensional subspace of linear forms of
Theorem~\ref{t:frommukai}, and let $G$ be its  $6$-dimensional quotient:
$$0 \rightarrow F \rightarrow \wedge^2 V \rightarrow G \rightarrow
0.$$
The space $F$ cuts out  $\proj G \subset \proj \wedge^2 V$, in which $C$ is
nondegenerately embedded.  Thus $\proj G$ is identified with 
$\proj H^0(C, \cK_C)$, and thus up to a scalar multiple $G$ is identified with
$H^0(C, \cK_C)$.

\point The vector bundle $E$ has a useful description in terms of the
$g^1_4$'s of $C$.

\tpoint{Theorem} {\em  If $\cL$ is a line bundle giving a $g^1_4$, and $\cM =
\cK _C\otimes \cL^{\vee}$ is its Serre-dual $g^2_6$, then $E$ may
be written in an exact sequence 
\begin{equation}
\label{eq:LEM}
\xymatrix{ 0 \ar[r] & \cL \ar[r]^\be & E \ar[r] & \cM \ar[r] & 0.}
\end{equation}
This is  the unique (up to scaling) nontrivial extension of $\cM$ by $\cL$.
}\label{t:frommukai2}\ravispace

The genus $8$ analog of Theorem~\ref{t:frommukai2} is shown in
\cite[Lem.~3.6]{mukai}, and as described in \cite[\S 5]{mukai}, the
analogous argument holds in genus $6$.

\point We will need a refined ``family'' version of
Theorem~\ref{t:frommukai}, which is undoubtedly known to experts,
although we were unable to find a reference.  This describes $\cm_6^M$
as a quotient $(Y \setminus \De) / GL(5)$.  We now define the terms in
this quotient.

Let $V$ be a dimension $5$ vector space.
Consider the Grassmannian $G(2,V)$, which has itsPl\"ucker embedding in
$\proj \wedge^2 V$.  The $4$-dimensional vector space  of linear sections of $G(2,V)$ (under
the Pl\"ucker embedding) are thus parametrized by $G(4, \wedge^2 V)$.

The quadratic forms on such complete intersections (again, in the
Pl\"ucker variables) are parametrized by a rank $16$ vector bundle  on
$G(4, \wedge^2 V)$: there is the $\binom {5+2} 2 = 21$-dimensional
vector space of quadratic forms on the $\proj^5$ cut out in $\proj
(\wedge^2 V)$ by the $4$ linear forms; but a $5$-dimensional vector
space of them lie in the ideal generated by the Pl\"ucker quadrics.

Denote this rank $16$ vector bundle over $G(4, \wedge^2 V)$ by $Y$.
Note that $\dim Y = 40$.  Then $Y$ parametrizes complete intersections
of $G(2,V)$ with four linear forms and a quadratic form.  The locus in
$Y$ where such complete intersections are singular is a divisor $\De$.
Then $Y \setminus \De$ corresponds to the smooth complete intersections.  It
is straightforward (see for example \cite[\S 5]{mukai}) that such
complete intersections are canonically embedded (irreducible) genus
$6$ curves.  The group $GL(V)$ acts on $Y \setminus \De$ through its
action on $G(2,V)$.

\tpoint{Theorem}  
{\em The natural map $\phi: (Y \setminus \De) / GL(V) \rightarrow
  \cm_6$
is an open embedding (of Deligne-Mumford stacks) whose image is
precisely $\cm_6^M$.}\label{t:mukai}\ravispace

For this reason, we call $\cm_6^M$ the {\em Mukai locus}.  Denote the
universal curve over $\cm_6^M$ by $\cC^M \rightarrow \cm_6^M$.  If
$[C] \in \cm_6^M$, we say that $C$ is a {\em Mukai-general} curve (of
genus $6$).

\bpf Note first that the image of $\phi$ is precisely $\cm_6^M$, by
Mukai's work.  By the
Enriques-Petri theorem (see \cite[Main Thm.]{saint-donat}), these are the only
canonical curves cut out by quadrics. Conversely, if a curve $C$ is
Mukai-general, then it is in the image of $\phi$ by
\cite[Thm.~5.2]{mukai}.  To find an inverse image of $C$ in $(Y
\setminus \De)$, we take the unique stable rank $2$ vector bundle $E$
with $\det E = K_C$ and $h^0(C,E) = 5$, which happens to be globally
generated (Theorem~\ref{t:frommukai}), and take the induced map $C
\rightarrow G(2,5)$.

We next claim that the map $\phi$ is representable, i.e., gives an
isomorphism of isotropy groups.  The previous paragraph establishes $Y
\setminus \De$ as representing the moduli space of triples of a curve
in the Mukai locus $C$, a stable rank $2$ vector bundle $E$ with
$\det(E) \cong \cK_C$ and $h^0(C, E)=5$, and a framing of
$h^0(C,E)$.  The quotient by $GL(5)$ forgets the framing.
We thus must show that for each curve $C$ in the Mukai locus, the
vector bundle $E$ is unique up to scaling.  This is a slight
strengthening of Mukai's result \cite[Thm.~5.1(1)]{mukai}, which shows
that $E$ is unique up to isomorphism.

We filter $E$ as in \eqref{eq:LEM}.
Hence given any automorphism $\al: E \rightarrow E$, we have a diagram
$$
\xymatrix{
0 \ar[r] & \cL  \ar[r]^\be &  E \ar[d]^\al  \ar[r] & \cM \ar[r] & 0 \\
0 \ar[r] & \cL \ar[r]^\be &  E  \ar[r] & \cM \ar[r] & 0. \\
}
$$
  We need the further fact that the  inclusions $\cL \hookrightarrow E$
are unique up to scaling.  For this we need that $\Hom(\cL, \cM) =
0$.
But $\Hom(\cL, \cM) = H^0(C, \cM \otimes \cL^{\vee})$ is $0$,
as
$\cM \otimes \cL^{\vee}$ is a line bundle of degree $5-3=2$, and $C$ is not
hyperelliptic.
Thus after rescaling $\al$, we can extend the diagram as follows:
$$
\xymatrix{
0 \ar[r] & \cL \ar[d]^=  \ar[r]^\be &  E \ar[d]^\al  \ar[r] & \cM \ar[r] & 0 \\
0 \ar[r] & \cL \ar[r]^\be &  E  \ar[r] & \cM \ar[r] & 0 \\
}
$$
Finally, noting that there is only one nontrivial extension of $\cM$
by $\cL$ 
(up to scaling of $\cL$ and $\cM$ and
the extension class in $\Ext^1(\cM, \cL)$) 
yielding a vector bundle with $E$ with $5$ sections
(\cite[Lem.~3.6]{mukai}, translated to the genus $6$ setting as in
\cite[\S 5]{mukai}), after rescaling $\cM$ in the bottom row of the
diagram above, we may write
$$
\xymatrix{
0 \ar[r] & \cL \ar[d]^=  \ar[r]^\be &  E \ar[d]^\al  \ar[r] & \cM \ar[d]^= \ar[r] & 0 \\
0 \ar[r] & \cL \ar[r]^\be &  E  \ar[r] & \cM \ar[r] & 0 \\
}
$$
from which we see that $E$ is unique up to scaling, as desired.

We have thus established that $\phi$ is representable. 
Finally, by Zariski's Main Theorem (which applies for Deligne-Mumford
stacks, see \cite[Thm.~C.1]{fantechi}), as $\cm_6^M$ is
normal,  the map $\phi$ is an isomorphism.
   \epf

\epoint{Remarks}
(i) Just as Mukai's paper \cite{mukai} treats genus $6$ curves by
precisely the same method as genus $8$ curves, the proof of
Theorem~\ref{t:mukai}  applies with the obvious changes to the genus
$8$ case.

(ii) For our applications, we need less than that $\phi$ is
representable; we need only that $\phi$ is a relative Deligne-Mumford
stack, as we are showing results about Chow groups with
$\Q$-coefficients.  However, it seems more elegant to show this
better result, rather than describing a technical work-around.

\bpoint{The Mukai vector bundle $\cV$ on $\cm_6^M$}
This argument shows the existence of a rank $5$ vector bundle $\cV$ on
$\cm_6^M$ (and a rank $2$ vector bundle $\cE$ on the universal curve $\cC^M$ over
$\cm_6^M$), relativizing Mukai's construction.  For obvious reasons we
call $\cV$ the {\em Mukai
bundle}.  It is disappointing
that our argument uses Zariski's Main Theorem (and thus in
particular the lucky fact that the moduli space of curves is normal);
there should be a direct argument.  But we were unable to directly  ``relativize''
Mukai's argument.

By the construction of $\cm_6^M$ of Theorem~\ref{t:mukai},
we have a tautological section to $G(4, \Sym^2 \cV) \rightarrow
\cm_6^M$, from which we have a tautological exact sequence of locally
free sheaves on $\cm_6^M$:
\begin{equation}\label{eq:FVE}
0 \rightarrow \cF \rightarrow \Sym^2 \cV \rightarrow \E' \rightarrow 0
\end{equation}
where $\rank \cF = 4$ and $\rank \E' = 6$.      
We use the notation $\E'$ because we may interpret it as a
``twist'' of the Hodge bundle $\E$:

\tpoint{Proposition} \quad \newline 
{\em 
(a) We have 
$\E' \cong \E \otimes \cL'$ for some invertible sheaf $\cL'$ on
$\cM_6^M$.  \newline
(b) The Chern classes of $\E'$ are tautological.}\label{p:Eprime}

\bpf For a point $[C] \in \cm_6^M$ of the Mukai locus, the fibers of $\cF$
and $\E'$ are the vector spaces $F$ and $G$ of \S \ref{s:frommukai}.
In particular, $\proj \E'$ is canonically isomorphic to $\proj \E$,
so that the following diagram commutes:
$$
\xymatrix{\cC^M \ar[r]  \ar[dr] &  \proj \E' \ar[d] \\
& \proj \E
}
$$
where the vertical map $\proj \E' \rightarrow \proj \E$ is the (canonical) isomorphism, the diagonal
map  $\cC^M \rightarrow \proj \E$ is the canonical map of $\cC^M$ into the projectivization of its
relative canonical bundle, and the horizontal map $\cC^M \rightarrow \proj \E'$ is the
map corresponding to the pointwise construction of \S
\ref{s:frommukai}.    This implies part (a).

 As $\Pic_{\Q} \cM_6^M \equiv \Pic_{\Q} \cM_6$ (they
differ by a codimension $2$ subset, by Theorem~\ref{t:6stra}), and 
$\Pic_{\Q} \cM_6$ is generated by $\la$ (Proposition~\ref{p:picm6}), we have part (b).
\epf

\section{The Mukai-general locus:  $A^*(\cm^M_6)$ is tautological}
\label{s:mukai}In this section 
we  complete the proof of the Main Theorem~\ref{t:main} by proving:

\tpoint{Theorem}
{\em $A^*(\cm^M_6)$ is tautological.}\label{t:muktau}\ravispace

For convenience, let $v_i := c_i(\cV)$ (the Chern classes of the Mukai
bundle), and $f_i := c_i(\cF_)$ (the Chern classes of the tautological
rank $4$ vector bundle $\cF$ on $G(4, \Sym^2 \cV)$).

\tpoint{Proposition}  {\em $A^*(\cm_6^M)$ is generated by $f_1$,
  \dots, $f_4$, $v_1$, \dots, $v_5$.}\label{p:fv}

\bpf
By Vistoli's
Theorem~\ref{t:vistoli}, the Chow ring of $\cm_6^M$ is generated by the
Chern classes of $\cV$ and the Chow groups of $Y \setminus \De$.

Now $Y$ is a projective bundle over $G(4,\Sym^2 \cV)$, and $\De$ is a
divisor of nonzero relative degree over $G(4, \Sym^2 \cV)$ (for each
point of $G(4, \Sym^2 \cV)$, we can find a quadric that restricts to a
singular curve on the surface corresponding to that point of $G(4,
\Sym^2 \cV)$), so the Chern classes of $Y \setminus \De$ are generated
by the Chow groups of $G(4, \Sym^2 \cV)$, which are in turn generated
as an algebra by the Chern classes of the tautological bundle.     
\epf

Thus to prove Theorem~\ref{t:muktau}, it suffices to prove that $f_1$,
\dots, $f_4$, $v_1$, \dots, $v_5$ are tautological.

By
exact sequence \eqref{eq:FVE},  the $f_i$ lie in $R^*(v_1, \dots, v_5)$, so
it suffices to prove that $v_1$, \dots, $v_5$ are tautological. 

Let $\pi: \cC \rightarrow \cm_6^M$ be the universal curve over the
Mukai-general locus.  Then on $\cm_6^M$, we have an exact sequence 
$$
0 \rightarrow \cG \rightarrow \Sym^2 \pi_* \omega \rightarrow \pi_*
\omega^{\otimes 2} \rightarrow 0.$$
Here $\omega$ is the relative dualizing sheaf as usual, and 
$\cG$ is a rank $6$ vector bundle, which over each point $[C] \in
\cm_6$ corresponds to the $6$-dimensional vector space of quadrics
cutting out the canonical curve $C \subset \proj H^0(C,
\omega_C)$.

Just as in \S \ref{s:g5}, the Chern classes of $\Sym^2 \pi_* \omega = \Sym^2 \E$ are 
tautological, and the Chern classes of $\pi_* \omega^{\otimes 2}$ are
tautological using Grothendieck-Riemann-Roch.  Thus the
Chern classes of $\cG$ are also tautological.

Note that $\cG$ has a canonical rank $5$ subbundle corresponding to
the $5$ quadrics cutting out the Pl\"ucker embedding of $G(2, V)$.  In
other words, we have an exact sequence
$$
0 \rightarrow \cV' \rightarrow \cG \rightarrow \cL \rightarrow 0
$$
where $\cL$ is an invertible sheaf.  Now $\cL$ has tautological Chern
class (by Proposition~\ref{p:picm6}), so $\cV'$ has tautological Chern
classes.

\tpoint{Lemma}  {\em Suppose $V$ is a dimension $5$ vector space.   
We have the Pl\"ucker embedding $G(2,V) \hookrightarrow \proj \wedge^2
V$, cut out by the 5 Pl\"ucker quadrics.  The quadrics are, as a
$GL(V)$-representation, $\wedge^4 V = V^\vee \otimes \det V$.}

\bpf This is a standard argument.  (For example, let $e_1$, \dots,
$e_5$ be a basis of $V$.  We may explicitly write down the
Pl\"ucker quadrics in terms of the Pl\"ucker variables $e_{ij} := e_i
\wedge e_j$.  Character theory then allows us to determine
the $GL(V)$-representation on these quadrics.)  \epf

By our construction of $\cm_6^M$ as a $GL(V)$-quotient, we have 
$\cV' = \cV^{\vee} \otimes \det \cV$.  As $\det \cV$ has tautological
Chern class (Proposition~\ref{p:picm6}), $\cV^\vee$ has tautological
Chern classes, and hence $\cV$ does too.

\end{document}